\documentclass[11pt]{article}
\usepackage{a4wide}

\usepackage{epsfig}
\usepackage{amsmath}
\usepackage{amssymb} 
\usepackage{makeidx}
\usepackage{graphics}
\usepackage{color}

\begin{document}
 
\def\RR {\mathbb{R}}
\def\NN {\mathbb{N}}
\def\TT {\mathbb{T}}  
\def\TU {\TT_{\rm U}} 
\def\MM {\mathbb{M}}  
\def\CC {\mathbb{C}}   
\def\II {\mathbb{I}}
\def\PP {\mathcal{P}}
\def\SS {\mathcal{S}}
\def\VV {\mathcal{V}}
\def\HH {\mathcal{H}}
\def\UU {\mathcal{U}}
\def\KK {\mathcal{K}}
\def\opspan {\mbox{\rm span}}
\newtheorem{Th}{Theorem}[section] 
\newtheorem{Le}[Th]{Lemma}
\newtheorem{Co}[Th]{Corollary} 
\newtheorem{Pro}[Th]{Proposition}
\newtheorem{Def}[Th]{Definition} 
\newtheorem{rem}[Th]{Remark}
\newcommand{\hdrie}{\hspace{3mm}}
\newcommand{\und}{\hdrie\mbox{\rm and }\hdrie}
\newcommand{\half}{\frac{1}{2}}
\newcommand{\be}{\begin{equation}}
\newcommand{\ee}{\end{equation}}
\newcommand{\diag}{\mbox{\rm diag}}
\newcommand{\sth}{\hdrie | \hdrie}
\newcommand{\trace}{{\rm trace}}
\newcommand{\ol}{\overline}
\newcommand{\Mnm}{\mathbb{M}^{n\times m}(\mathbb{\CC})}
\newcommand{\Mmm}{\mathbb{M}^{m\times m}(\mathbb{\CC})}
\newcommand{\Mpp}{\mathbb{M}^{p\times p}(\mathbb{\CC})}
\newcommand{\Mnp}{\mathbb{M}^{n\times p}(\mathbb{\CC})}
\newcommand{\Mkk}{\mathbb{M}^{k\times k}(\mathbb{\CC})}
\newcommand{\Mnn}{\mathbb{M}^{n\times n}(\mathbb{\CC})}
\newcommand{\with}{\hdrie\mbox{\rm with }\hdrie}
\newcommand{\Ret}{\Theta}
\newcommand{\Imt}{\Theta^\perp}
\newcommand{\eqv}{\hspace{1mm}\Leftrightarrow\hspace{1mm}}
\newcommand{\st}{\sqrt{2}} 
\newcommand{\hst}{\half\sqrt{2}} 

\title{Normalitity preserving perturbations and augmentations and their effect on the eigenvalues}

\author{Ricardo Reis da Silva and Jan H. Brandts}
\date{\today}

\maketitle

\begin{abstract} We revisit the normality preserving augmentation of normal matrices studied by Ikramov and Elsner in 1998 and complement their results by showing how the eigenvalues of the original matrix are perturbed by the augmentation. Moreover, we construct all augmentations that result in normal matrices with eigenvalues on a quadratic curve in the complex plane, using the stratification of normal matrices presented by Huhtanen in 2001. To make this construction feasible, but also for its own sake, we study normality preserving normal perturbations of normal matrices. For $2\times 2$ and for rank-$1$ matrices, the analysis is complete. For higher rank, all essentially Hermitian normality perturbations are described. In all cases, the effect of the perturbation on the eigenvalues of the original matrix is given. The paper is concluded with a number of explicit examples that illustrate the results and constructions.
\end{abstract}

{\bf Keywords:} structured perturbation theory, normal matrix, essentially Hermitian matrix, Toeplitz decomposition, normal augmentation problem.

\section{Introduction}
A well-known sufficient condition for the sum of normal matrices $A$ and $E$ to be normal is, that $A$ and $E$ commute. This condition, however, is not necessary. If $A$ and $E$ are both Hermitian, the sum is also trivially normal. In both the above situations, a well-developed perturbation theory for the eigenvalues of $A+E$ in terms of $A$ and its perturbation $E$ is available. See for instance \cite{StSu} and the references therein. For other {\em normality preserving normal perturbations}, one could of course apply one of the many equivalent characterizations \cite{IkEl,GrJoSaWo} of the normality of $A+E$, but other than that there seems to be much less literature. Although Wielandt \cite{Wie} already studied the location of eigenvalues of sums of normal matrices in 1955, he did not require the sum to be normal.

In this paper, we study the perturbation of normal matrices by {\em essentially Hermitian matrices} \cite{BeKoPr,Dru, Fre}. These are matrices $E$ that can be written as $E=\beta H+\alpha I$, where $\alpha,\beta\in\CC$, $H$ is Hermitian, and $I$ is the identity matrix. Since all $2\times 2$ normal matrices and all rank-$1$ normal matrices are essentially Hermitian, this will completely characterize the $2\times 2$ and the rank-$1$ case, but will also give insight in the $k\times k$ and rank-$(k-1)$ case for $k\geq 3$. Moreover, we will relate the eigenvalues of $A+E$ to those of $A$.

Essentially Hermitian matrices have been frequently discussed, both in the numerical and the core linear algebra setting. See for instance the introductory section of \cite{Fre}, where the discretization of the Helmholz equation $\Delta u + k^2u = f$ where $k=k_1+k_2i\in\CC$ is the wave number and $u$ the unknown eigenmode of the differential operator, results in essentially Hermitian matrices. Faber and Manteuffel show in \cite{FaMa} that for linear systems with essentially Hermitian system matrix, there exists a variant of the Conjugate Gradient method that still relies on a three term recursion. In the context of the Arnoldi method applied to normal matrices, Huckle \cite{Huc} proves that irreducible normal tridiagonal matrices are essentially Hermitian. In fact, the upper Hessenberg matrices generated by the Arnoldi method is tridiagonal only for essentially Hermitian matrices. The behavior of Arnoldi's method is in that case similar as the Lanczos method. Generally, however, the Arnoldi method for normal matrices does not respect the structure. Huhtanen in \cite{Huh1} raises this issue. He shows that for every normal matrix $A$ and almost all unimodular $z\in\CC$ the skew-Hermitian part $S$ of $zA$ is a polynomial $p$ of degree at most $n-1$ in the Hermitian part $H$ of $zA$. The essentially Hermitian matrices are exactly the ones for which the degree of $p$ is at most one, resulting in collinear eigenvalues. Since the $p$ can be retrieved in a modest number of arithmetic operations as a by-product of the Lanczos algorithm applied to $H$, this led to the development \cite{Huh1,Huh2} of efficient structure preserving algorithms for eigenvalue problems and linear systems involving normal matrices.

Essentially Hermitian matrices also appear in the context of core linear algebra problems. For instance, in \cite{Ikr}, Ikramov investigates which matrices $R$ can be the upper triangular part of a normal matrix $A$. Although the full problem remains open, a solution exists if all diagonal entries of $R$ lie on the same line in the complex plane, say
\be \ell: \RR\rightarrow\CC: \hdrie \rho \mapsto y + \theta \rho, \hdrie y\in\CC \und \theta\in\CC, \hdrie |\theta|=1. \ee
Writing $R=D+U$ where $D$ is the diagonal of $R$, we see that $\ol{\theta}(D+U-yI)$ is upper triangular with real diagonal, and consequently, $H=\ol{\theta}(D+U-yI) + \theta U^*$ is Hermitian. Therefore,
\be A = \theta H + y I = R+\theta^2 U^*\ee
is essentially Hermitian, hence normal, with upper triangular part equal to $R$. Together with Elsner \cite{IkEl2}, Ikramov then considered the {\em normality preserving augmentation problem}: when does the addition rows and columns to a given normal matrix $A$ yield a normal matrix $A_+$? Essentially Hermitian matrices appear in this setting as special cases, though the problem in its full generality remains unsolved. We will take the augmentation problem as a starting point for this paper. We add to their results an analysis of the eigenvalues of the augmented matrix $A_+$ in terms of $A$. Moreover, we investigate how Huhtanen's stratification of the normal matrices may be of help in providing additional solutions. Our analysis involves normality preserving normal perturbations as mentioned in the beginning of this introduction, a topic that we will study in more detail as well. Also in this context we present results on the eigenvalues of the perturned matrix. 

\subsection{Outline}
In Section \ref{sect1} we define counterparts of the Hermitian and skew-Hermitian parts of a matrix, based on unimodular complex numbers $\theta$ and $\theta i$, that for $\theta=1$ reduce to the standard Hermitian and skew-Hermitian part. This will somewhat facilitate the discussion on essentially Hermitian matrices of the form $\theta H+\alpha I$. We also outline the description of normal matrices as given in \cite{Huh1} and pay attention to commuting matrices and their eigenvalues and eigenvectors. In Section \ref{sect2}, we revisit the normality preserving augmentation problem by Ikramov and Elsner \cite{IkEl} and add to their analysis results on the eigenvalues of the resulting augmented matrix. In Section \ref{sect3} we use similar ingredients to study the normality preserving normal perturbation problem. This will be of use in Section \ref{sect5}, where we explicitly construct augmentations of normal matrices that have all eigenvalues on a quadratic curve in $\CC$.

\subsection{Notation}
For quick reference, we list here some notations that are frequently used in this paper.
\begin{itemize}
\item $\Mnm$ is the set of all $n\times m$ matrices with complex entries,
\item $I$, or $I_n\in\Mnn$ is the identity matrix with columns $e_1,\dots,e_n$,
\item $\sigma(A)$ the set of eigenvalues of $A$, its {\em spectrum},
\item $\TT$ is the circle group of unimodular numbers,
\item $\TU$ is the subset of $\TT$ in the closed upper half plane, minus $-1$,
\item $\Ret(A)$ and $\Imt(A)$ are $\theta$-(skew)-Hermitian parts of $A$ (see Section \ref{sect1.1}),
\item $[A,B] = AB-BA$ is the matrix commutator.
\item $\PP^k(\RR)$ is the space of real-valued polynomials on $\RR$ of degree less than or equal to $k$
\end{itemize}
The zero matrix we simply denote by $0$ regardless of its dimensions, and within matrices we will often use empty space to indicate a zero block.
  
\section{Preliminaries}\label{sect1}
In this section we define a family of non-standard decompositions of complex numbers and matrices and prove some of their properties. We also recall a result on the eigenvalues and eigenvectors of commuting normal matrices.

\subsection{The $\theta$-Toeplitz decomposition and some relevant properties}\label{sect1.1}
A complex number $z\in\CC$ is often split up as $z = \Re(z) + \Im(z)$, where $2\Re(z) = z + \ol{z}$ and $2\Im(z) = z - \ol{z}$. This interprets the complex plane as a two-dimensional real vector space with as basis the numbers $1$ and $i$. In this paper, it will be convenient to decompose $z$ differently. For this, we introduce a family of decompositions parametrized in $\theta\in\TT\subset\CC$, the circle group of unimodular numbers. For a given $\theta\in\TT$ we let
\be\label{een} z = \Ret(z) + \Imt(z), \hdrie\mbox{\rm where }\hdrie \Ret(z) = \Re(\ol{\theta}z)\theta \und  \Imt(z) = \Im(\ol{\theta}z)\theta.\ee 
Moreover, apart from the standard Toeplitz or Cartesian decomposition of a square matrix $A$ into its Hermitian and skew-Hermitian parts,
\be A = \HH(A) + \SS(A),\hdrie\mbox{\rm where }\hdrie \HH(A) = \half(A + A^*) \und  \SS(A) = \half(A - A^*),\ee 
in accordance with (\ref{een}), we consider the family of matrix decompositions
\be\label{theta} A = \Ret(A)+\Imt(A),\hdrie\mbox{\rm where }\hdrie \Ret(A) = \HH(\ol{\theta}A)\theta \und  \Imt(A) = \SS(\ol{\theta}A)\theta.\ee
We will call this decomposition of $A$ its {\em $\theta$-Toeplitz} decomposition, $\Ret(A)$ the {\em $\theta$-Hermitian part} and $\Imt(A)$ the {\em $\theta$-skew-Hermitian} part of $A$ and subsequently, $A$ is {\em $\theta$-Hermitian} if $A=\Ret(A)$ and {\em $\theta$-skew-Hermitian} if $A=\Imt(A)$. Note that $\theta$-skew Hermitian matrices are $i\theta$-Hermitian. We now generalize a well-known result to the  $\theta$-Toeplitz decomposition of normal matrices.

\begin{Le}\label{lemma1} Let $\theta\in\TT$ be arbitrarily given. For any normal matrix $A$ we have that
\be\label{twob} Av = \lambda v \hdrie \Leftrightarrow \hdrie \Ret(A)v = \Ret(\lambda)v \und \Imt(A)v = \Imt(\lambda)v, \ee 
where $\Ret$ and $\Imt$ and their relation to $\theta$ are defined in (\ref{theta}).
\end{Le}
{\bf Proof. } Let $Au=\lambda u$ for some $u\not=0$. Since $A$ is normal, there exists a unitary matrix $U$ with $u$ as first column and $U^*AU=\Lambda$ diagonal. But then $U^*A^*U=\Lambda^*$, showing that $A^*u=\ol{\lambda}u$. Since the argument can be repeated with $A^*$ instead of $A$, this yields that $A$ and $A^*$ have the same eigenvectors and that the corresponding eigenvalues are each other's complex conjugates. Therefore,
\be 2\Ret(A)v = (\ol{\theta}A+\theta A^*)\theta v = (\ol{\theta}\lambda+\theta\ol{\lambda})\theta v = 2\Ret(\lambda)v\ee
and, similarly, also $\Imt(A)v=\Imt(\lambda)v$. The reverse implication is trivial.\hfill $\Box$

\begin{Co}\label{cor1} Let $\lambda_1$ be an eigenvalue of a normal matrix $A$. For given $\theta\in\TT$, consider the line $\ell\subset\CC$ through $\lambda_1$ defined by
\be \ell: \RR\rightarrow \CC: \rho\mapsto \lambda_1 + \rho\theta. \ee
Assume that $\lambda_1,\dots,\lambda_p$ are all eigenvalues of $A$ that lie on $\ell$. Then the eigenspace $\UU$ of the eigenvalue $\Imt(\lambda_1)$ of $\Imt(A)$ equals the invariant subspace of $A$ spanned by $u_1,\dots,u_p$. Restricted to $\UU$, the matrix $A-\lambda_1 I_p$ is $\theta$-Hermitian.
\end{Co} 
{\bf Proof.} Let $\lambda_a,\lambda_b$ be eigenvalues of $A$, then we have that
\be \Imt(\lambda_a) = \Imt(\lambda_b) \eqv \Im(\ol{\theta}\lambda_a) = \Im(\ol{\theta}\lambda_b) \eqv \ol{\theta}(\lambda_a-\lambda_b) = \rho\in\RR \eqv \lambda_a-\lambda_b = \theta \rho.\ee
Thus, $\Imt(A)u_j=\Imt(\lambda_1)u_j$ for all $j\in\{1,\dots,p\}$, and conversely, if $\Imt(A)u=\Imt(\lambda_1)u$ then $u$ is a linear combination of $u_1,\dots,u_p$. Writing $U_p$ for the matrix with columns $u_1,\dots,u_p$ we moreover find that
\be AU_p = U_p\Lambda_p \hdrie\mbox{\rm with }\hdrie \Lambda_p = \theta R+ \lambda_1I_p\ee
where $\Lambda_p$ is the $p\times p$ diagonal matrix whose eigenvalues are $\lambda_1,\dots,\lambda_p$ and 
$R$ is real diagonal. Thus, $\Imt(U_p^*AU_p-\lambda_1 I_p)=0$, proving the last statement.  \hfill $\Box$

\subsection{Normal matrices with all eigenvalues on polynomial curves}\label{sect2.2}
If $A$ is $\theta$-Hermitian then $A$ is normal because $\ol{\theta}A$ is Hermitian, and by the spectral theorem for Hermitian matrices, all eigenvalues of $A$ lie on the line $\ell:\RR\rightarrow\CC: \rho\mapsto \rho\theta$. In the literature, for instance in \cite{BeKoPr,Dru, Fre}, $A$ is called {\em essentially Hermitian} if there exists an $\alpha\in\CC$ such that $A-\alpha I$ is $\theta$-Hermitian for some $\theta\in\TT$. Clearly, the spectrum of an essentially Hermitian matrix lies on an affine line shifted over $\alpha\in\CC$. Conversely, if a normal matrix has all its eigenvalues on a line $\ell\subset\CC$, it is essentially Hermitian. This includes {\em all} normal $2\times 2$ matrices and {\em all} normal rank one perturbations of $\alpha I$ for $\alpha\in\CC$. Larger and higher rank normal matrices have their eigenvalues on a polynomial curve $\mathcal{C}\subset\CC$ of higher degree. 
 
\subsubsection{Polynomial curves of degree $k\geq 2$}
Each normal $A\in\Mnn$ has its eigenvalues on polynomial curve $\mathcal{C}\subset\CC$ of degree $k\leq n-1$. This can be explained as follows \cite{Huh1}. First, fix $\theta\in\TT$ such that for each pair $\lambda_p,\lambda_q$ of eigenvalues of $A$
\be\label{str1} \lambda_p\not=\lambda_q \Rightarrow \Ret(\lambda_p)\not=\Ret(\lambda_q).  \ee
Note that there exist at most $n(n+1)$ values of $\theta$ for which this cannot be realized, corresponding to the at most $\half n(n+1)$ straight lines through each pair of distinct eigenvalues of $A$. Once (\ref{str1}) is satisfied, the points
\be\label{interp} (\ol{\theta}\Ret(\lambda_1),\ol{i\theta}\Imt(\lambda_1)), \dots, (\ol{\theta}\Ret(\lambda_n),\ol{i\theta}\Imt(\lambda_n))\in \RR\times \RR \ee
form a feasible set of points in $\RR\times\RR$ through which a Lagrange interpolation polynomial $\pi\in\PP^{n-1}(\RR)$ can be constructed that satisfies
\be \Imt(\lambda_j) = i\cdot\theta\pi(\ol{\theta}\Ret(\lambda_j)) \hdrie \mbox{\rm for all $j\in\{1,\dots,n\}$}.\ee
Notice that $\theta\pi(\ol{\theta}\Ret(\lambda_j))$ is $\theta$-Hermitian. The following Lemma summarizes the consequences.  

\begin{Le} Let $A$ be normal and $\theta$ such that (\ref{str1}) is satisfied. Then there is a $\pi\in\PP^{n-1}(\RR)$ such that $\Imt(A) = i\cdot\theta\pi(\ol{\theta}\Ret(A))$, and thus the $\theta$-Toeplitz decomposition of $A$ can be written as
\be A = \Ret(A)+i\cdot\theta\pi(\ol{\theta}\Ret(A)),\ee
or, in terms of $N=\ol{\theta}A$ and its classical Toeplitz decomposition, $N = \HH(N) + i\cdot\pi (\HH(N))$. The eigenvalues of $A$ lie on the image $\mathcal{C}$ of the function
\be c:\RR\rightarrow \CC: \hdrie \rho \mapsto \theta\rho + i\cdot\theta\pi(\rho). \ee
\end{Le}
{\bf Proof. } Apply the spectral theorem for normal and Hermitian matrices. \hfill $\Box$
 
\begin{figure}[ht]
\centering
\includegraphics[width=11cm]{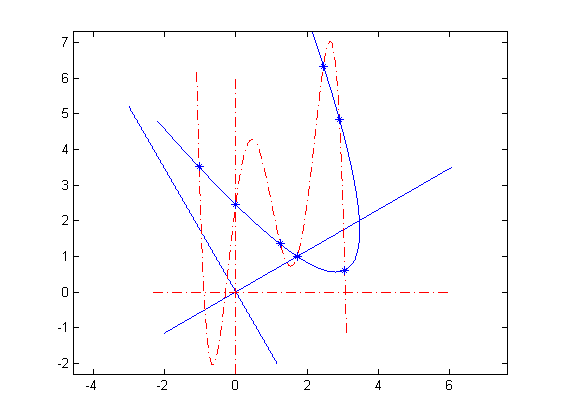}
\caption{Interpolating polynomials of different degree for different values of $\theta$. The seven asterisks represent the eigenvalues.}\label{figure1}
\end{figure}

\begin{rem}{\rm Note that $A$ is essentially Hermitian if and only if there exists a $\theta\in\TT$ such that the interpolating polynomial  $\pi\in\PP^1(\RR)$. In fact, $\pi$ may then even be in $\PP^0(I)$.}
\end{rem}

\begin{rem}{\rm If for some normal matrix $A$ the degree of the interpolation polynomial equals two for some value of $\theta$, it may well be of degree $n-1$ for almost all other values of $\theta$, since then the eigenvalues lie on a rotated parabola, as depicted in Figure \ref{figure1}. There does not seem to be an easy way to determine $\theta$ for which the polynomial degree is minimal, although for each given $\theta$ the polynomial can be computed exactly in a finite number of arithmetic operations without knowing the eigenvalues.}
\end{rem}
 
\subsubsection{Computing the polynomial $\pi$ for given $\theta\in\TT$}\label{sect1.2.2}
For almost any fixed value of $\theta$, the interpolation polynomial $\pi$ belonging to a normal matrix $A$ can be computed, without knowing the eigenvalues of $A$, in a finite number of arithmetic operations, providing us with a curve $\mathcal{C}\subset\CC$ on which all eigenvalues of $A$ lie. Indeed, if the degree of $\pi$ equals $k$ then $A-\Ret(A)$ is a linear combination of 
\be i\theta I,\,i\theta\Ret(A),\,i(\theta\Ret(A))^2,\,\dots,\,i(\theta\Ret(A))^{k}.\ee
Making the combination explicit is equivalent to finding the coefficients of $\pi$. To obtain $\pi$ in practice, notice that for any $v\in\CC^n$,
\be \pi(\ol{\theta}\Ret(A))v \in \KK^k(\ol{\theta}\Ret(A),v) = \opspan\{v,\ol{\theta}\Ret(A)v, \dots , (\ol{\theta}\Ret(A))^kv\},  \ee
the Krylov subspace for $\ol{\theta}\Ret(A)$ and $v$. Since $\ol{\theta}\Ret(A)$ is Hermitian, an orthonormal basis for $\KK^k(\ol{\theta}\Ret(A),v)$ can be constructed using a three-term recursion, and solving the linear system can be done cheaply. These, and other considerations, led Huhtanen to the development of efficient structure preserving eigensolvers \cite{Huh1} and linear system solvers \cite{Huh2} for problems involving normal matrices.
 
\subsection{Commuting normal matrices}
Because we will need to draw conclusions about the eigenvalues and eigenvectors of commuting normal matrices, here we recall a well-known result and give its complete proof.

\begin{Le}\label{lemma3} Normal matrices $A$ and $E$ commute if and only if they are simultaneously unitarily diagonalizable.
\end{Le}
{\bf Proof.}  If both $\Lambda=W^*AW$ and $\Delta=W^*EW$ are diagonal for a unitary matrix $W$, then clearly $[A,E] = W[\Lambda,\Delta]W^* = 0$. Conversely, assume that $[A,E]=0$. Let $U$ be a unitary matrix such that $\Delta =U^*EU$ is diagonal with multiple eigenvalues being neighbors on the diagonal of $\Delta$. Thus
\[ \Delta = \left[ \begin{array}{ccc} \delta_1 I_{m_1} & & \\ & \ddots & \\ & & \delta_\ell I_{m_\ell}\end{array}\right],\]
where $m_j$ denotes the multiplicity of $\delta_j$. Let $S=U^*AU$ and write $s_{pq}$ for its entries. Then equating the entries of $S\Delta$ and $\Delta S$ in view of the relation
\[ [S,\Delta] = U^*[A,E]U = 0, \] 
shows that $s_{pq}=0$ whenever $\delta_p\not=\delta_q$. Thus, $S$ is block-diagonal with respective blocks $S_1,\dots,S_\ell$ of sizes $m_1,\dots,m_\ell$. For each $j\in\{1,\dots,\ell\}$, $S_j$ is normal. Let $Q_j$ be unitary and such that $\Lambda_j = Q_j^*S_jQ_j$ is diagonal. Then with
\be Q = \left[ \begin{array}{ccc} Q_1 & & \\ & \ddots & \\ & & Q_\ell\end{array}\right] \und \Lambda = \left[ \begin{array}{ccc} \Lambda_1 & & \\ & \ddots & \\ & & \Lambda_\ell\end{array}\right] \und W=UQ  \ee
we find that $\Delta = Q^*\Delta Q = W^*EW$ and $\Lambda=Q^*SQ=W^*AW$.\hfill$\Box$

\begin{rem}\label{remmult}{\rm In case all eigenvalues of $E$ are distinct, then $S$ itself is already diagonal and the proof is finished. In case $E$ has an eigenvalue, say $\delta_1$, of multiplicity $m_1>1$, then there is freedom in the choice for the first $m_1$ columns $u_1,\dots,u_{m_1}$ of $U$ that correspond to $\delta_1$. Even though each choice diagonalizes $E$, not each choice diagonalizes $A$ as well. This is expressed as $S$ having a diagonal block $S_1$ of size $m_1$. Writing $U_1$ for the matrix with columns $u_1,\dots,u_{m_1}$ we have that 
\be AU_1=U_1S_1, \ee
hence the column span of $U_1$ is an invariant subspace $\VV$ of $A$. The matrix $Q_1$ determines, through the transformation $W_1 = U_1Q_1$, an orthonormal basis for $\VV$ of eigenvectors of $A$. If $S_1$ has multiple eigenvalues, again there may be much freedom in the choice for $Q_1$.}
\end{rem}

\section{Normality preserving augmentation}\label{sect2} 
In this section we revisit, from an alternative point of view, a problem studied by Ikramov and Elsner in [11]. It concerns the augmentation by a number $m$ of rows and columns of a normal matrix in such a way, that normality is preserved. Our analysis differs from the one in \cite{IkEl2}, and we add details on the eigendata of $A_+$ in terms of those of $A$. 
 
\bigskip\noindent
{\bf Normality preserving augmentation. } {\em Let $A\in\Mnn$ be normal. Characterize all $n\times m$ matrices $V,W\in\Mnm$, and all $Y\in \Mmm$ such that $A_+$, where
\be\label{augment} A_+ = \left[ \begin{array}{lr} A & V \\ W^* & Y \end{array}\right],  \ee
is normal, too. In other words, characterize all normality preserving augmentations of $A$}.
 
\begin{rem} {\rm Note that Hermicity, $\theta$-Hermicity, and essentially Hermicity preser\-ving augmentation problems are all trivial, because each of these properties is inherited by principle submatrices. For unitary matrices, this does not hold. However, if $A$ and $A_+$ are both unitary, their rows and columns all have length one. Thus $V=W=0$ and also $Y$ is unitary. This trivially solves the unitarity preserving augmentation problem. When we consider normality preserving augmentation, matters become less trivial.}
\end{rem}

\subsection{The normality preserving augmentation problem for $m=1$}\label{sect2.1}
First consider the case $m=1$, and write $v,w,y$ instead of $V,W,Y$. It is easily verified that $[A_+,A_+^*]=0$ if and only if
\be\label{five}  ww^* = vv^*, \hdrie w^*w = v^*v, \und A^*v+ yw = A w + \ol{y} v. \ee
The two leftmost relations hold if and only if $v=\phi w$ for some $\phi\in\TT$. The rightmost relation may add further restrictions on $v,w$ and $\phi$. Before studying these, note however that $\phi$ is the square of a unique $\theta\in\TU\subset\TT$, where 
\be \TU = \{\, \tau \in \TT \sth \arg(\tau) \in [0,\pi) \,\}. \ee
This yields a reformulation of $v=\phi w$ that better reveals the underlying structure, which is
\be\label{seven} v = \phi w \hdrie\Leftrightarrow\hdrie v = \theta^2 w \hdrie\Leftrightarrow\hdrie \ol{\theta}v = \theta w = u\hdrie\Leftrightarrow\hdrie v = \theta u \und w = \ol{\theta} u \ee
for some $u\in\CC^n$. Further restrictions on $u$ and $\theta\in\TU$ follow from substituting $v=\theta u$ and $w=\ol{\theta}u$ into the rightmost equation in (\ref{five}), which after some rearrangements results in the requirement that
\be\label{ess} \Imt(A)u = \Imt(y)u.\ee 
Because eigenpairs of $\Imt(A)$ were already characterized in Corollary \ref{cor1}, we have solved the augmentation problem for $m=1$. The below theorem summarizes this solution constructively in terms of the eigendata of $A$. Note that this result was proved already in \cite{IkEl2}, though in a different manner.

\begin{Th}\label{th1} Let $y\in\CC$, and $\ell:\RR\rightarrow\CC: \rho\mapsto y+\rho\theta$ a line in $\CC$ through $y$ with slope $\theta\in\TU$. Let $\lambda_1,\dots\lambda_p$ be the eigenvalues of $A$ that lie on $\ell$. Then, the matrix 
\be A_+ = \left[ \begin{array}{rl} A & v \\ w^* & y\end{array}\right] \ee
is normal if and only if $v=\theta u$ and $w=\ol{\theta}u$, where $u$ is a linear combination of eigenvectors corresponding to $\lambda_1,\dots\lambda_p$, with the convention that $u=0$ if $p=0$.
\end{Th} 
{\bf Proof. } Corollary \ref{cor1} shows the relation between eigendata of $\Imt(A)$ and $A$, and together with the derivation in this Section this proves the statement.\hfill $\Box$

\subsection{The eigenvalues of the augmented matrix}
Next, we augment the analysis in \cite{IkEl2} with a study of the eigenvalues of $A_+$ in relation to those of $A$. Let $\Lambda_p\in\Mpp$ be the diagonal matrix whose eigenvalues are the $p$ eigenvalues of $A$ that lie on $\ell:\RR\rightarrow\CC: \rho\mapsto y + \rho\theta$. Then, as already mentioned in the proof of Corollary \ref{cor1},
\be \Lambda_p = \theta R+ yI_p\ee
for some real diagonal matrix $R$. Let $U\in\Mnn$ be any unitary matrix whose last $p$ columns are eigenvectors of $A$ belonging to $\lambda_1,\dots,\lambda_p$ and let $U_p$ contain those last $p$ columns of $U$. Then, assuming that $A$ and $A_+$ are normal, Theorem \ref{th1} shows, with $r=U_p^*u$, that
\be\label{ci1} \left[\begin{array}{rl} U & \\ & 1\end{array}\right]^*\left[ \begin{array}{cl} A & v \\ w^* & y\end{array}\right]\left[\begin{array}{cc} U & \\ & 1\end{array}\right] =  \left[\begin{array}{clc} B & &  \\ & \Lambda_p & \theta r \\  & \theta r^* & y\end{array}\right].\ee
Moreover,
\be\label{ci2} \left[\begin{array}{lc} \Lambda_p  & \theta r \\ \theta r^* & y\end{array}\right] = 
\theta R_+ + y I_{p+1},\hdrie \mbox{\rm where } \hdrie R_+=\left[\begin{array}{lc} R & r \\ r^* & 0\end{array}\right]. \ee
The above observations reveal some additional features of the solution of the augmentation problem, that we will formulate as another theorem.

\begin{Th}\label{Th2} The only normality preserving $1$-augmentations of $A$ are the ones that, on an orthonormal basis of eigenvectors of $A$, augment a $p\times p$ essentially Hermitian submatrix of $A$. Hence, $n-p$ eigenvalues of $A$ are also eigenvalues of $A_+$. The $p+1$ eigenvalues of $A_+$ that remain lie on the same line as, and are interlaced by, the remaining $p$ eigenvalues of $A$.
\end{Th}
{\bf Proof. } The block form in (\ref{ci1}) shows that the eigenvalues of $B$ are eigenvalues of both $A$ and $A_+$, whereas (\ref{ci2}) shows that to locate the remaining $p+1$ eigenvalue of $A_+$ one only needs to observe that the real eigenvalues of $R$ interlace those of $R_+$.\hfill $\Box$

\begin{rem}{\rm Note that the case $p=0$, covered by Theorem \ref{th1}, is also included in the above analysis if one is willing to interpret {\em on the same line as the remaining $p$ eigenvalues of $A$} as any line in $\CC$. This just reflects that the additional eigenvalue $y\in\CC$ of $A_+$ can lie anywhere.}
\end{rem}

For an illustration of the constructions of Theorems \ref{th1} and \ref{Th2}, see Section \ref{ill1}. There, we will augment a given $3\times 3$ matrix in two different ways and compute the eigenvalues of the augmented matrix.
 
\subsection{Normal matrices with normal principal submatrices}
By applying the procedure for $m=1$ several times consecutively, also $m$-augmentations with $m>1$ can be constructed. In particular, {\em all} the normal matrices having the property that all their principal submatrices are normal, can be constructed. Since generally, normal matrices do not have normal principal submatrices, this shows that the $m$-augmentation for $m>1$ has not yet been completely solved. In Section \ref{sect5} we will study the principal submatrices of normal matrices from the point of view of Section \ref{sect2.2}. This will also give more insight into the $m$-aug\-men\-ta\-tion problem for $m=2$, that already proved to be very complicated in \cite{IkEl2}. In particular, we will give a procedure to augment $A$ that does not reduce to $m$-fold application of the $1$-augmentation. Before that, we investigate normality preserving normal perturbations. Apart from being of interest on its own, it will be needed in Section \ref{sect5}.

\section{Normality preserving normal perturbations}\label{sect3}
In this section we consider a question related to the augmentation problem, and we will study it using the same techniques. In particular, we will study normality preserving $\theta$-Hermitian perturbations, that play a role in the augmentation problem in Section \ref{sect3.2}.
 
\bigskip\noindent
{\bf Normality preserving normal perturbation. } {\em Let $A\in\Mnn$ be normal. Characterize all normal $E$ such that $A_+=A+E$ is normal. In other words, characterize the normality preserving normal perturbations $E$ of $A$}.
   
\bigskip
Since $A=\Ret(A)+\Imt(A)$, and both $\Ret(A)$ and $\Imt(A)$ are normal, this shows that the sum of normal matrices can be literally {\em any} matrix and thus that the above problem is non-trivial. 
 
\subsection{Generalities}
To start, we formulate a multi-functional lemma that summarizes the technicalities of writing out commutators of linear combinations of matrices.

\begin{Le}\label{lem1} Let $A,E\in\Mnn$ and $\gamma,\mu\in\CC$. Then with $\theta=\ol{\gamma}\mu/|\gamma\mu|$,
\be\label{a1} [\gamma A+\mu E, (\gamma A+\mu E)^*] = |\gamma|^2[A,A^*]+ 2\gamma\ol{\mu}\cdot\Ret([A,E^*])+|\mu|^2[E,E^*],\ee
and
\be\label{et} \Ret([A,E^*]) = [\Ret(A),\Imt(E^*)] + [\Imt(A),\Ret(E^*)]. \ee
Therefore, if $A$ and $E$ are normal, then $\gamma A+\mu E$ is normal if and only if 
\be \label{twtw} \Ret([A,E^*]) = 0,  \ee
or in other words, if and only if $[A,E^*]$ is $\ol{\gamma}\mu$-skew-Hermitian.
\end{Le}
{\bf Proof.} Straight-forward manipulations with the commutator give the statements. \hfill $\Box$

\begin{Co}\label{cor2}{\sl Let $A,E\in\Mnn$ be normal. Then $A+E$ is normal if and only if $\gamma A + \mu E +\alpha I$ is normal for all $\gamma,\mu,\alpha\in\CC$ with the restriction that $\gamma\ol{\mu}\in\RR$.} 
\end{Co} 

\begin{Co} If $A,E\in\Mnn$ are normal then $[A,E^*]=0 \Leftrightarrow [A,E]=0$. Thus if either term vanishes, both $A+E^*$ and $A+E$ are normal.
\end{Co} 
{\bf Proof. } As was shown in the proof of Lemma \ref{lemma1}, $E$ and $E^*$ have the same eigenvectors. Thus, $A$ and $E^*$ are simultaneously unitarily diagonalizable if and only if $A$ and $E$ are. Lemma \ref{lemma3} now proves that $[A,E^*]=0 \Leftrightarrow [A,E]=0$, and Lemma \ref{lem1} proves the conclusion. \hfill $\Box$ 
 
\begin{Co}\label{Co4.4} The sum of  $\gamma A$ and $\mu E$ with $\gamma,\mu\in \CC$ and $A$ and $E$ Hermitian is normal if and only if $\ol{\gamma}\mu\in\RR$ or $[A,E]=0$.
\end{Co} 
{\bf Proof. } The commutator of Hermitian matrices is always skew-Hermitian. Thus for $[A,E]$ to be $\theta$-Hermitian in (\ref{twtw}), $\ol{\gamma}\mu$ must be real, or $[A,E]$ should vanish. \hfill $\Box$
 
\subsection{Normality preserving normal rank one perturbations}
This section aims to show the similarities between the normality preserving normal rank one perturbation problem and the $m$-augmentation problem of Section \ref{sect2.1}. Indeed, let $E=vw^*$ with $v,w\in\CC^n$. Then $E$ is normal if and only if 
\be  \|w\|^2 vv^* = \|v\|^2 ww^*, \ee
and thus if and only if $v=zw$ for some $z\in\CC$. Write $z=\theta \rho$ with $\theta\in\TT$ and $\rho\in\RR, \rho\geq 0$. This shows that a rank-$1$ matrix $E$ is normal if and only if $E$ is $\theta$-Hermitian,  
\be E = \theta uu^*, \hdrie\theta\in\TT. \ee  
With $A\in\Mnn$ normal, we will look for the conditions on $u\in\CC^n$ and $\theta\in\TT$ such that $A+\theta uu^*$ is normal. Since $\theta uu^*$ is $\theta$-Hermitian, by (\ref{et}) in Lemma \ref{lem1},
\be\label{two} [(A+\theta uu^*),(A+\theta uu^*)^*] = 2\ol{\theta}[\Imt(A),uu^*], \ee 
thus, $A+\theta uu^*$ is normal if and only if $\Imt(A)$ and $uu^*$ commute. According to Lemma \ref{lemma3}, this is if true and only if they are simultaneously unitarily diagonalizable. For this, it is necessary and sufficient that $u$ is an eigenvector of $\Imt(A)$. As in Theorem \ref{th1} we will formulate this result constructively in terms of the eigendata of $A$.

\begin{Th}\label{Th1} Let $y\in\CC$, and $\ell:\RR\rightarrow\CC: \rho\mapsto y + \rho\theta$ a line in $\CC$ through $y$ with slope $\theta\in\TT$. Let $\lambda_1,\dots\lambda_p$ be the eigenvalues of $A$ that lie on $\ell$. Then, the matrix 
\be A_+ = A + \theta uu^* \ee 
is normal if and only if $u$ is a linear combination of eigenvectors corresponding to $\lambda_1,\dots\lambda_p$, with the convention that $u=0$ if $p=0$.
\end{Th} 
{\bf Proof. } Corollary \ref{cor1} shows the relation between eigendata of $\Imt(A)$ and $A$, and together with the derivation above this proves the statement.\hfill $\Box$

\begin{rem}\label{iterate}{\rm An interesting consequence of adding the normality preserving normal rank one perturbation $E=\theta uu^*$ is that
\be\label{eq10} \Imt(A+\theta uu^*) = \Imt(A) +\Imt(\theta uu^*) = \Imt(A),\ee
because $\theta uu^*$ is $\theta$-Hermitian. Thus, the conditions under which adding {\em another} $\theta$-Hermitian rank one perturbation $F=\theta ww^*$ to $A+E$ lead to a normal $A+E+F$ are {\em identical} to the conditions just described for $E$. We will get back to this observation in Section \ref{sect3.3}. }
\end{rem}

\subsection{The eigenvalues of the perturbed matrix}
To study the eigenvalues of $A_+$ in relation to those of $A$, let $\Lambda_p\in\Mpp$ be the diagonal matrix whose eigenvalues are the $p$ eigenvalues of $A$ that lie on $\ell:\RR\rightarrow\CC: \rho\mapsto y+\rho\theta$. Then
\be \Lambda_p = \theta R+ yI_p\ee
for some real diagonal matrix $R$. Let $U\in\Mnn$ be any unitary matrix whose first $p$ columns are eigenvectors of $A$ belonging to $\lambda_1,\dots,\lambda_p$. Then, assuming that $A$ and $A_+$ are normal, Theorem \ref{Th1} shows that
\be U^*A_+U = U^*(A+\theta uu^*)U =  \left[\begin{array}{cc} \Lambda_p &  \\ & B \end{array}\right] + \left[\begin{array}{cc} \theta rr^* & \\  & \end{array}\right].\ee
because $u$ is a linear combination of the first $p$ columns of $U$. This leads to the following theorem, in which we summarize the above analysis.

\begin{Th}\label{Th4.7} The only normality preserving normal rank one perturbations of $A$ are the ones that, on an orthonormal basis of eigenvectors of $A$, are $\theta$-Hermitian rank one perturbations of a $p\times p$ $\theta$-Hermitian submatrix. Hence, $n-p$ eigenvalues of $A$ are also eigenvalues of $A_+$. The remaining $p$ eigenvalues of $A_+$ are the eigenvalues of
\be\label{eigs} \Lambda_p + \theta rr^* = \theta (R+rr^*) + yI_p, \ee
and these interlace the $p$ eigenvalues of $A$ on $\ell$ with the additional $p+1$-st point $+\infty \theta$.
\end{Th}
{\bf Proof. } The eigenvalues of $A_+$ are the eigenvalues of $B$ together with the eigenvalues of the matrix in (\ref{eigs}). Obviously, all eigenvalues of $B$ are eigenvalues of $A$ as well. Since $rr^*$ is a positive semi-definite rank one perturbation of $R$, the eigenvalues $\rho_1\leq \cdots\leq \rho_p$ of $R+rr^*$ and the eigenvalues $r_1\leq \cdots r_p$ of $R$ satisfy
\be r_1 \leq \rho_1 \leq r_2 \leq \cdots \leq \rho_{p-1} \leq r_p \leq \rho_p \ee
as a result of Weyl's Theorem \cite{StSu}. Multiplying by $\theta$ and shifting over $y$ yields the proof. \hfill $\Box$ 

\begin{rem}{\rm Note that if $p=1$, only one eigenvalue is perturbed, and $\rho_1 = r_1 + \|r\|^2$. In terms of the original perturbation $E=\theta uu^*$ this becomes $\tilde{\lambda}=\lambda + \|u\|^2$, where $\lambda$ is the eigenvalue of $A$ belonging to the eigenvector $u$.}
\end{rem}
As a consequence of the following theorem, it is possible to indicate where the eigenvalues of the family of matrices $A+tE$ are located. This can only be done for {\em normal} normality preserving perturbations.

\begin{Th}\label{Th5}{\sl Let $A,B\in\Mnn$ be normal. Consider the line $\ell$ through $A$ and $B$,
\be \ell: \RR\rightarrow \Mnn: \hdrie t \mapsto tA+(1-t)B.\ee
If $E=B-A$ is normal, all matrices on $\ell$ are normal; if $E$ is not normal, $A$ and $B$ are the only normal matrices on $\ell$}.
\end{Th}
{\bf Proof. } Observe that $\ell(t)=A+(1-t)E$. If $E$ is normal, then Corollary \ref{cor2} shows that all matrices $\gamma A+\mu E$ with $\gamma\ol{\mu}\in\RR$ are normal, which includes the line $\ell$. Assume now that $E$ is not normal. Because $A$ and $B=A+E$ both are normal, (\ref{a1}) in Lemma \ref{lem1} gives that 
\be (1-t)\left(2\HH([A,E^*])+(1-t)[E,E^*]\right) = 0\with [E,E^*]\not=0. \ee
The solution $t=1$ confirms the normality of $A$, and the linear matrix equation
\be 2\HH([A,E^*])+(1-t)[E,E^*] = 0 \with [E,E^*]\not=0\ee
allows at most one solution in $t$ which, by assumption, is $t=0$.\hfill$\Box$\\[2mm]
Thus, any line in $\Mnn$ parametrized by a real variable that does not lie entirely in the set of normal matrices, contains at most two normal matrices.
 
\begin{rem}\label{remeigs}{\rm Lemma \ref{lemma1} shows that if $E\in\Mnn$ is such that $A+E$ is normal, then 
\be \sigma(A+E) \subset \sigma(\Ret(A)+\Ret(E)) \times \sigma(\Imt(A)+\Imt(E)), \ee
and perturbation theory for $\theta$-Hermitian matrices can be used to derive statements about the eigenvalues of $A+E$. According to Theorem \ref{Th5}, if $E$ {\em itself} is normal too, this relation is valid {\em continuously} in $t$ along the line $A+tE$:
\be \sigma(A+tE) \subset \sigma(\Ret(A)+t\Ret(E)) \times \sigma(\Imt(A)+t\Imt(E)). \ee
For non-normal $E$ this would generally not be true, as will be illustrated in Section \ref{ill2}.}
\end{rem}

\begin{Co}\label{Co4.10} As a result of Theorem \ref{Th5} and Remark \ref{remeigs}, the perturbed eigenvalues of 
\be t\mapsto A + t\theta uu^*, \hdrie 0\leq t \leq 1,\ee
seen as functions of $t$, are line segments that all lie on the same line with slope $\theta$.
\end{Co}
 
\subsection{$\theta$-Hermitian rank-$k$ perturbations of normal matrices}\label{sect3.3}
In Section \ref{sect3.2} we will encounter normality preserving $\theta$-skew-Hermitian perturbations. In this section we will fully characterize those. For this, consider for given $k$ with $1\leq k\leq n$ the $\theta$-Hermitian rank-$k$ matrix
\be\label{eq11} E = \theta H, \hdrie \mbox{\rm where } \hdrie \theta\in\TU \und H=H^*.\ee 
Let $A$ be normal. Then, since $\Imt(E^*)=0$, Lemma \ref{lem1} shows that $A+E$ is normal if and only if
\be\label{eq48} [\Imt(A),H] = 0.\ee
By Lemma \ref{lemma3} this is equivalent to $\Imt(A)$ and $H$ being simultaneously diagonalizable by a unitary transformation $U$. Thus, $H$ needs to be of the form
\be\label{delta}  H = U\Delta U^* \ee
where $\Delta\in\mathbb{M}^{k\times k}(\RR)$ is diagonal with diagonal entries $\delta_1,\dots,\delta_k$ and the columns $u_1,\dots,u_k$ of $U$ are orthonormal eigenvectors of $\Imt(A)$. But then, writing
\be E = E_1+\cdots +E_k, \hdrie\mbox{\rm where for all $j\in\{1,\dots,k\}$, }\hdrie E_j = \theta\delta_ju_ju_j^*, \ee
the observation in Remark \ref{iterate} reveals that perturbing $A$ by $E$ is equivalent to perturbing $A$ consecutively by the rank one matrices $E_1,\dots,E_k$. The above analysis is summarized in the following theorem. An illustration of this theorem is provided in Section \ref{ill3}.
 
\begin{Th}\label{Th4} Let $E=\theta H$ be a $\theta$-Hermitian rank $k$ perturbation of a normal matrix $A$. Then $E$ is normality preserving  if and only if $E$ can be decomposed as 
\be E = E_1 + \cdots + E_k, \ee 
where $E_1,\dots, E_k$ are all normality preserving $\theta$-Hermitian rank one perturbations of $A$. In fact, for each permutation $\sigma$ of $\{1,\dots,k\}$ and each $m\in\{1,\dots,k\}$, the partial sum
\be A + \sum_{j=1}^m E_{\sigma(j)} \ee
is then normal, too.
\end{Th}

\begin{rem}{\rm In accordance with Remark \ref{remmult}, if $\Delta$ in (\ref{delta}) has multiple eigenvalues, there exist non-diagonal unitary matrices $Q$ such that $\Delta=Q\Delta Q^*$. As a result, $H$ can be written as $Z\Delta Z^*$ where the orthonormal columns of $Z$ span an invariant subspace of $\Imt(A)$. This implicitly writes the perturbation $\theta H$ as a sum of rank-$1$ normal perturbations that do not necessarily preserve normality. This aspect is also illustrated in Section \ref{ill3}.}
\end{rem}

\begin{Th}\label{Th4.8} The only normality preserving $\theta$-Hermitian rank-$k$ perturbations of $A$ are the ones that, on an orthonormal basis of eigenvectors of $A$, are $\theta$-Hermitian perturbations 
of $\theta$-Hermitian submatrices of size $s_1\times s_1, \dots, s_m\times s_m$ of ranks $k_1,\dots,k_m$, where $k_1+\cdots+k_m=k$. As a result of this perturbation, at most $s_1+\cdots+s_m$ eigenvalues of $A$ are perturbed, which are located on at most $m$ distinct parallel lines $\ell_1,\dots,\ell_m$, defined by $\theta\in\TT$ and $y_1,\dots,y_m\in\CC$ as
\be \ell_j:\RR\rightarrow\CC: \hdrie \rho \mapsto y_j + \theta\rho. \ee
Moreover, the eigenvalues of $A+tE, t\in[0,1]$ connect the perturbed eigenvalues of $A$ with those of $A$ by line segments that lie on $\ell_1,\dots,\ell_m$.
\end{Th} 
{\bf Proof. } Write $\theta H = \theta(\delta_1u_1u_1^*+\dots+\delta_k u_ku_k^*)$, where $u_1,\dots,u_k$ are eigenvectors of $\Imt(A)$, and repeatedly apply Theorem \ref{Th4.7}. The statement about the eigenvalues of $A+tE$ follows from Corollary \ref{Co4.10}.\hfill $\Box$

\begin{figure}[ht]
\centering 
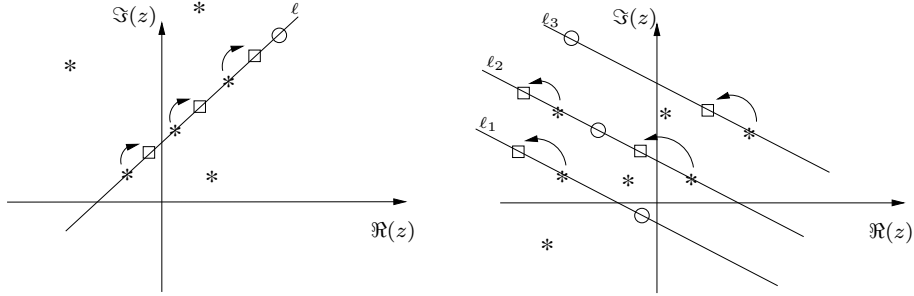
\caption{Eigenvalue perturbation by a rank-$1$ matrix (left) and a rank-$k$ matrix (right).}\label{figure2}
\end{figure}

\bigskip
For a qualitative illustration of the effect on the eigenvalues due to a rank-$1$ $\theta$-Hermitian perturbation and a rank-$k$ $\theta$-Hermitian perturbation of a normal matrix, see Figure \ref{figure2}. The asterisks in the pictures are eigenvalues of $A$, the circles represent different choices for $y$, and the boxes are the perturbed eigenvalues.

\begin{rem}\label{rem4.9}{\rm If $H$ in (\ref{eq11}) is semi-definite, then the eigenvalues of $A+tE$ all move in the same direction over the parallel lines $\ell_1,\dots,\ell_m$ from Theorem \ref{Th4.8}.}
\end{rem} 

\subsection{Normality preserving normal perturbations of normal matrices}
The above analysis of $\theta$-Hermitian normality preserving perturbations also gives sufficient conditions for when normal perturbations $E$ of the form $E=\theta_1H_1+\theta_2H_2$ with $H_1,H_2$ Hermitian and $\theta_1,\theta_2\in\TT$ are normality preserving in case $\theta_2$ does not equal $\pm\theta_1$. Notice that in order for $E$ to be normal itself, $[H_1,H_2]=0$ by Corollary \ref{Co4.4}. Obviously, $E$ is in general not $\theta$-Hermitian for some value of $\theta$. Nevertheless, the following holds.

\begin{Th}\label{Thx} Let $A$ be normal, $H_1,H_2$ Hermitian with $[H_1,H_2]=0$, and $\theta_1,\theta_2\in\TT$ with $\theta_1\ol{\theta_2}\not\in\RR$. Then $E=E_1+E_2$ is a normality preserving normal perturbation of $A$ if $E_1$ and $E_2$ both are normality preserving perturbations of $A$.
\end{Th} 
{\bf Proof. } Corollary \ref{Co4.4} covers the normality of $E$. Furthermore, assuming that $A+E_1$ is normal, (\ref{eq48}) gives that $A+E_1+E_2$ is normal if and only if $[\Imt(A+E_1),H_2]=0$, where $\Imt$ indicates the $\theta_2$-skew Hermitian part. But
\be [\Imt(A+E_1),H_2] = [\Imt(A),H_2] + [\Imt(\theta_1H_1),H_2] = [\Imt(A),H_2],\ee
because $[H_1,H_2]=0$, proving the statement. \hfill $\Box$ 

\begin{rem}{\rm A similar result holds for normal perturbations $E=E_1+\cdots+E_k$ where $E_j=\theta_jH_j$ with $H_j$ Hermitian and $\theta_j\in\TT$ for each $j\in\{1,\dots,k\}$. Moreover, each normal perturbation $E$ can be written in this form in several different ways.}
\end{rem} 
 
We are now ready to return to the augmentation problem of Section \ref{sect3} and to study augmentations of normal matrices $A$ that result in an augmented matrix $A_+$ whose eigenvalues all lie on the graph of a quadratic polynomial.

\section{Further augmentations}\label{sect5}
We now return to the $m$-augmentation problem of Section \ref{sect3}. We concentrate on the case $m>1$ and on augmentations $A_+$ that have some additional structure in the spirit of Section \ref{sect2.2}.  For instance, if all eigenvalues of $A_+$ lie on a line, then $A_+$ is essentially Hermitian, a property that is inherited by principal submatrices. In that case it is clear which matrices $A$ can be augmented into $A_+$. The next simplest case, if simple at all, is the case where all eigenvalues of $A_+$ lie on  the image of a {\em quadratic} function in a rotated complex plane.

\subsection{Augmentations $A_+$ with all eigenvalues on a quadratic curve}\label{sect3.2}
Assume that for $\theta\in\TT$ there exists a quadratic $\pi\in\PP^2(\RR)$ such that 
\be\label{aug2}  A_+ = \Ret(A_+) + i\cdot \theta \pi\left(\ol{\theta}\Ret(A_+)\right), \hdrie\mbox{\rm where } \hdrie  A_+ = \left[ \begin{array}{lr} A & V \\ W^* & Y \end{array}\right]. \ee
Then $A_+$ is normal with all its eigenvalues on the rotated parabola $\mathcal{C}\subset\CC$ defined as the image of $q$, where
\be\label{curve} q: \RR\rightarrow\CC: \hdrie \rho \mapsto \theta\rho + i\cdot\theta \pi(\rho) \ee
and 
\be   \pi(x)=r_0+r_1x+r_2x^2 \hdrie\mbox{\rm with } \hdrie r_0,r_1,r_2\in\RR.\ee

\begin{rem}\label{remr}{\rm Throughout this section we assume, without loss of generality, that $r_2>0$. The case $r_2=0$, as argued above, is trivial and concerns essentially Hermitian matrices, whereas the case $r_2<0$ can be avoided by replacing $\theta$ by $-\theta$, which is nothing else than a trivial change of coordinates that transforms the polynomial $\pi$ into $-\pi$. }
\end{rem}
The curve $\mathcal{C}$ now divides the complex plane $\CC$ in three disjoint parts
\be\label{partc} \CC = \mathcal{C}_+\cup \mathcal{C}\cup\mathcal{C}_-,\ee
where $\mathcal{C}_+$ is the open part of $\CC$ that lies on the one side of $\mathcal{C}$ that is convex. 
 
\subsubsection{The principal submatrices of $A_+$ and their eigenvalues} \label{sect5.1.1}
For convenience, write
\be X=\ol{\theta}\Ret(A),\hdrie  M=\ol{\theta}\Ret(Y) \und 2Z = \ol{\theta}V+\theta W, \ee
then
\be \label{blocks} \ol{\theta}\Ret(A_+) = \left[\begin{array}{lr} X & Z \\ Z^* & M\end{array}\right],\und 
 (\ol{\theta}\Ret(A_+))^2 = \left[\begin{array}{rr} X^2+ZZ^* & XZ+ZM  \\ Z^*X+MZ^* & M^2+Z^*Z\end{array}\right].\ee   
Thus, explicitly evaluating $\pi$ at $A_+$ using the block forms in (\ref{blocks}), and comparing the result with the block form of $A_+$ displayed in (\ref{aug2}) yields that
\be\label{at} A = \Ret(A)+i\cdot\theta\pi(\ol{\theta}\Ret(A)) + i\cdot\theta r_2 ZZ^* \ee
and 
\be Y = \Ret(Y)+i\cdot\theta\pi(\ol{\theta}\Ret(Y)) + i\cdot\theta r_2 Z^*Z. \ee
Results that follow will sometimes be stated for $A$ only, even though similar statements obviously hold for $Y$. The first proposition simply translates (\ref{at}) in words.

\begin{Pro} The $n\times n$ leading principal submatrix $A$ of $A_+$ is a $\theta$-skew-Hermitian rank-$k$ perturbation of a normal matrix that has all its eigenvalues on $\mathcal{C}$, and $k\leq \min(m,n)$.
\end{Pro} 

\begin{Le}\label{lem4.10} If $A$ in (\ref{at}) is normal, then $\sigma(A)\subset\mathcal{C}\cup\mathcal{C}_+$.
\end{Le}  
{\bf Proof. } If $A$ is normal, then obviously $i\cdot\theta r_2 ZZ^*$ is a {\em normality preserving} $\theta$-skew-Hermitian perturbation of the normal matrix $\Ret(A)+i\cdot\theta\pi(\ol{\theta}\Ret(A))$ that has all its eigenvalues on $\mathcal{C}$. By Theorem \ref{Th4.8}, each eigenvalue $\lambda\in\mathcal{C}$ that is perturbed, moves along the line $\ell:\RR\rightarrow \CC: \rho\mapsto \lambda +i\cdot\theta\rho$. Note that $\ell$ is vertical in the $\theta$-rotated complex plane. By Remark \ref{rem4.9}, since $ZZ^*$ is positive semi-definite, the direction is the same for each perturbed eigenvalue, and is determined by the sign of $r_2$. In Remark \ref{remr} we assumed that $r_2>0$, and thus the direction is directed into $\CC_+$ defined in (\ref{partc}). \hfill $\Box$

\begin{rem}{\rm Note that a {\em multiple} eigenvalue $\lambda$ of $\Ret(A)+i\cdot\theta\pi(\ol{\theta}\Ret(A))$, located on $\mathcal{C}$, may be perturbed by $i\cdot\theta r_2 ZZ^*$ into several {\em distinct} eigenvalues of $A$. Those will all be located on $\ell:\RR\rightarrow \CC: \rho\mapsto \lambda +i\cdot\theta\rho$ with $\rho>0$}.
\end{rem}
 
\begin{Co}\label{Cor9} Assume that $A$ in (\ref{at}) is normal. Then $\sigma(A)\subset\mathcal{C}$ if and only $A_+$ is block diagonal with blocks $A$ and $Y$. 
\end{Co}
{\bf Proof. } If $Z\not=0$ then $\trace(ZZ^*)\not=0$ and at least one eigenvalue is perturbed. Lemma \ref{lem4.10} shows that a perturbed eigenvalue cannot stay on $\mathcal{C}$ and necessarily moves from $\mathcal{C}$ into $\mathcal{C}_+$. \hfill $\Box$

\subsubsection{Augmentations with eigenvalues on a quadratic curve}
The observations in Section \ref{sect5.1.1} can be reversed in the following sense. Given $A$, we choose a parabolic curve $\mathcal{C}$ and construct $Z\in\Mnm$ such that $i\cdot\theta r_2 ZZ^*$  perturbs the eigenvalues of $A$ onto $\mathcal{C}$. Then we use $\mathcal{C}$ to define the corresponding $m$-augmentation $A_+$ of $A$.

\begin{Co}\label{Co4.11} Necessary for a normal matrix $A$ to be $m$-augmentable into a normal matrix $A_+$ with all eigenvalues on a quadratic curve $\mathcal{C}$ is that $\sigma(A)\subset\mathcal{C}\cup\mathcal{C}_+$.
\end{Co}  
{\bf Proof. } This is just another corollary of Lemma \ref{lem4.10}.\hfill $\Box$

\bigskip
Clearly, for any given finite set of points in $\CC$, there are infinitely many candidates for such quadratic curves $\mathcal{C}$. It is the purpose of this section to show that {\rm each} of this candidates can be used, and to construct essentially all possible corresponding augmentations $A_+$.

\begin{Th}\label{Th7} Let $A\in\Mnn$ be normal, and let $\theta\in\TT$ and $\pi\in\PP^2(\RR)$ be such that $\sigma(A)\subset\mathcal{C}\cup\mathcal{C}_+$, where $\mathcal{C}$ is the graph of
\be \label{curve2}  q: \RR\rightarrow\CC: \hdrie \rho \mapsto \theta\rho + i\cdot\theta \pi(\rho). \ee
Then there exist $p$-augmentations $A_+$ of $A$ such that
\be \sigma(A_+)\subset\mathcal{C}, \ee
where $p$ is the number of eigenvalues of $A$ in $\mathcal{C}_+$. 
\end{Th}  
{\bf Proof. } Write $\Lambda_p\in\Mpp$ for the diagonal matrix with precisely the eigenvalues $\lambda_1,\dots,\lambda_p$ of $A$ that do not lie on $\mathcal{C}$ and let $U_p\in\Mnp$ have corresponding orthonormal eigenvectors $u_1,\dots,u_p$ as columns. Since $\sigma(\Lambda_p)\subset\mathcal{C}_+$, for each $j\in\{1,\dots,p\}$ there exists a positive real number $\xi_j$ such that  
\be \lambda_j-i\cdot\theta \xi_j \in\mathcal{C}.\ee
Write $\Xi\in\Mpp$ for the diagonal matrix with $\sigma(\Xi_p)=\{\sqrt{\xi_1},\dots,\sqrt{\xi_p}\}$ and set
\be Z = U_p\Xi_p.\ee 
By Lemma \ref{lemma1}, the columns of $U_p$ are also eigenvectors of $\Ret(A)$ and thus $i\cdot \theta ZZ^*$ is a $\theta$-skew-Hermitian normality preserving perturbation of $A$. Moreover,
\be u_j^*(A - i\cdot\theta ZZ^*)u_j =  \lambda-i\cdot\theta\xi \in\mathcal{C}. \ee
Because $A-i\cdot\theta ZZ^*$ is normal with all eigenvalues $\mathcal{C}$, the equality $\Ret(A-i\cdot\theta ZZ^*)=\Ret(A)$ leads to
\be A-i\cdot\theta ZZ^* = \ol{\theta}\Ret(A) + i\cdot\theta \pi(\ol{\theta}\Ret(A)). \ee
The assumption $r_2>0$, justified in Remark \ref{remr}, now gives that
\be A = \ol{\theta}\Ret(A) + i\cdot\theta \pi(\ol{\theta}\Ret(A)) + i\cdot\theta r_2\left(\frac{Z}{\sqrt{r_2}}\right)\left(\frac{Z}{\sqrt{r_2}}\right)^*, \ee
and according to (\ref{at}) this is precisely the $n\times n$ leading principal submatrix of $A_+$, where $A_+$ is defined as $A_+ = \ol{\theta}H + i\cdot\theta \pi(H)$, where
\be \label{ha} H = \left[\begin{array}{cc} \ol{\theta}\Ret(A) & \hat{Z} \\ \hat{Z}^* & M\end{array}\right], \hdrie\mbox{\rm with }\hdrie \hat{Z} = \frac{ZQ}{\sqrt{r_2}}\ee
and  $M,Q\in\Mmm$ are arbitrary Hermitian and unitary matrices. \hfill $\Box$

\begin{figure}[ht]
\centering
\includegraphics[width=11cm]{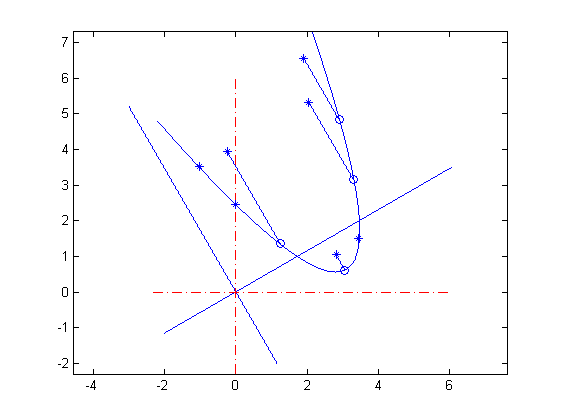}
\caption{Illustration of the construction in the proof of Theorem \ref{Th7}. Already $3$ of the $7$ eigenvalues of $A$ lie on the quadratic curve $\mathcal{C}$, and a rank-$4$ matrix $Z$ is needed to push the remaining $4$ eigenvalues onto $\mathcal{C}$.}\label{figure3}
\end{figure} 

\bigskip
For a given $n\times n$ normal matrix $A$, the typical situation is that after selecting $\theta$ suitably, at least $p$ of its eigenvalues do not lie on a quadratic curve, and a matrix $Z$ of rank $p$ is needed to push those outliers onto $\mathcal{C}$. In Figure 3, already $3$ of the $7$ eigenvalues of $A$, indicated by asterisks, lie on the quadratic curve $\mathcal{C}$, and a rank-$4$ matrix $Z$ is needed to push the remaining $4$ eigenvalues onto $\mathcal{C}$. After that, the augmented matrix can be formed.

\begin{rem}{\rm It is, of course, possible to move each eigenvalue of $A$ from $\mathcal{C}_+$ onto $\mathcal{C}$ as a result of an arbitrary amount of rank-$1$ perturbations. This would increase the number of columns of $Z$, and give $m$-augmentations of $A$ with $m>p$. However, this would not increase the rank of $Z$, and $ZZ^*$ would remain the same. Together with the analysis of Section \ref{sect3.3}, that shows which $\theta$-Hermitian perturbations are normality preserving, this shows that in essence, each $p$-augmentation $A_+$ with $Z$ of full rank, is of the form (\ref{ha}). In Section \ref{ill4} we give an explicit example of the construction in the proof of Theorem \ref{Th7}}.
\end{rem}

\subsubsection{Augmentations without computing eigenvalues}
So far, explicit knowledge about the eigenvalues and eigenvectors of $A$ was used to construct augmentations $A_+$. There are, however, cases in which it is sufficient to know the polynomial curve $\mathcal{E}$ on which the eigenvalues of $A$ lie. To see this, assume that $A$ is normal and
\be A = \Ret(A) + i\cdot\theta\pi(\ol{\theta}\Ret(A)), \hdrie \pi\in\PP^{2k}(\RR) \und \pi\not\in\PP^{2k-1}(\RR)\ee
for some integer $k\geq 1$. Since $\pi$ has even degree, there exist polynomials $p\in\PP^2(\RR)$ such that
\be \pi(x) - p(x) \geq 0 \hdrie \mbox{\rm for all $x\in\RR$}.\ee
This implies that the matrix $(\pi-p)(\ol{\theta}\Ret(A))$ is positive semi-definite, and hence it can be factorized as 
\be\label{factor} (\pi-p)(\ol{\theta}\Ret(A)) = ZZ^*, \ee
after which we have that
\be A = \Ret(A) + i\cdot\theta p(\ol{\theta}\Ret(A)) + i\cdot\theta ZZ^*.\ee
It is trivial that $i\cdot\theta ZZ^*$ is a normality preserving perturbation of $\Ret(A) + i\cdot\theta p(\ol{\theta}\Ret(A))$, and by choosing between $\theta$ and $-\theta$, as explained in Remark \ref{remr}, this leads to an $m$-augmentation of $A$, with generally $m=n-1$. Section \ref{sect1.2.2} explained how $\pi$ can be computed in a finite number of arithmetic operations, and the same is valid for the factorization (\ref{factor}). Of course, the problem of finding a minorizing polynomial $p\in\PP(\RR)$ may prove to be difficult in specific situations.

\subsection{Polynomial curves of higher degree}
If one tries to generalize the approach of Section \ref{sect3.2} to polynomial curves $\mathcal{C}$ of higher degree, the situation becomes rapidly more complex. As an illustration, consider the cubic case. The third power of the matrix $\ol{\theta}\Ret(A_+)$ in (\ref{blocks}) equals
\be \left[\begin{array}{cc} X^3 +XZZ^* + ZZ^*X + ZMZ^* & X^2Z+XZM+ZM^2+ZZ^*Z \\
Z^*X^2 + Z^*ZZ + MZ^*X + M^2Z & M^3 + Z^*ZM + MZ^*Z + Z^*XZ\end{array}\right], \ee
and thus, comparing the leading principal submatrices, 
\be A = \Ret(A)+i\cdot\theta\pi(\ol{\theta}A) + i\cdot\theta r_2ZZ^* + i\cdot\theta r_3(XZZ^* + ZZ^*X + ZMZ^*),\ee
where $r_3$ is the coefficient of $x^3$ of $\pi$. Thus, $A$ is a rank-$k$, with $k\leq 2m$, $\theta$-skew-Hermitian perturbation of the normal matrix $\Ret(A)+i\cdot\theta\pi(\ol{\theta}A)$. Of course, if $ZZ^*$ commutes with $\Ret(A)$, then it commutes with $X$, and this may help the analysis. However, it becomes much harder to control the perturbation in such a way, that $A$ will be augmented into a matrix $A_+$ with $\sigma(A_+)\subset\mathcal{C}$. Therefore, we will not pursue this idea in this paper.
 
\section{Illustrations}
In this section we present some illustrations of the main constructions and theorems of this paper. By making them explicit, we hope to create more insight in their structure.

\subsection{Illustrations belonging to Section \ref{sect2}}\label{ill1}
\begin{figure}[ht]
\centering
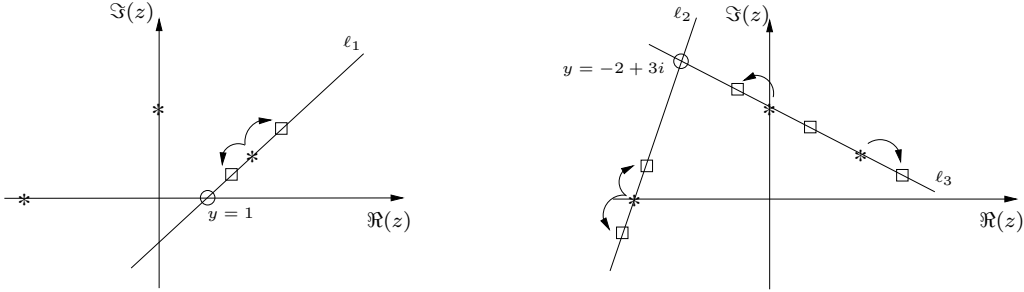
\caption{Left: the choice $y=1$ and the line $\ell_1$ trough the eigenvalue $2+i$. Right: the two lines $\ell_2$ and $\ell_3$ are depicted that correspond to the choice $y=-2+3i$.}\label{figure4}
\end{figure}

This example illustrates Theorems \ref{th1} and \ref{Th2}. Let $A$ be the matrix 
\be A = \left[\begin{array}{ccr} 2i &0 & 0\\ 0& 2+i & 0\\ 0& 0& -3\end{array}\right]. \ee
Take $y=1$ and choose $\ell_1$ to be the line through $y=1$ and the eigenvalue $2+i$ of $A$,
\be \ell_1: \RR\rightarrow \CC: \hdrie \rho \mapsto 1+\theta\rho, \hdrie\mbox{\rm with }\hdrie \theta = \frac{1+i}{\sqrt{2}}.\ee
Thus, in order for $A_+$ to be normal, $u$ must be a multiple of $e_2$, and $v=\theta u, w^*=\theta u^*$, which yields that for all $\mu\in\CC$,
\be A_+ = \left[\begin{array}{ccr|c} 2i & 0& 0&0\\ 0& 2+i &0 & \theta\mu \\0 &0 & -3 & 0\\ \hline 0& \theta\ol{\mu} &0 &  1\end{array}\right] \ee
is a normal augmentation of $A$. Moreover, for $y=1$ and this choice for $\ell_1$, these are all the normal augmentations of $A$. The eigenvalues of $A_+$ that are not eigenvalues of $A$ are the eigenvalues of
\be   \left[\begin{array}{cc} 2+i & \theta\mu\ \\ \theta\ol{\mu} & 1 \end{array}\right] =  \left[\begin{array}{cc} 1+\sqrt{2}\theta & \theta\mu\ \\ \theta\ol{\mu} & 1 \end{array}\right] = \theta \left[\begin{array}{cc} \sqrt{2} & \mu \\ \ol{\mu} & 0\end{array}\right] + I, \ee
and thus equal to
\be \lambda = \theta \left( \frac{\sqrt{2}}{2} \pm \sqrt{2}\mu^2 + 1\right) + 1, \ee
which lie on $\ell_1$ and have the eigenvalue $2+i$ that was perturbed as average, as is depicted in the left picture of Figure \ref{figure4}. The stars in both pictures are eigenvalues of $A$, the open circles represent different choices for $y$. Perturbed eigenvalues are indicated by squares.

A second option is to choose $y=-2+3i$ instead of $y=1$. This gives other possibilities to construct normality preserving augmentations. The first one is to choose $\ell_2$ through $y$ and $-3$, which shows that $u$ must be a multiple of $e_3$ and $v=\theta u, w^*=\theta u$ where $\theta = (1+3i)/\sqrt{10}$, which is a similar situation as for $y=1$. The second non-trivial option is to choose $\ell_3$ through $y$ and both $2i$ and $2+i$. Then with $\theta = (-2+i)/\sqrt{5}$, we may take $u$ as any linear combination of $e_1$ and $e_2$, showing that
\be A_+ = \left[\begin{array}{ccr|c} 2i & 0& 0& \theta \gamma \\ 0& 2+i & 0& \theta \mu \\ 0&0 & -3 & 0\\ \hline \theta \ol{\gamma} & \theta\ol{\mu} & 0& -2+3i \end{array}\right] \ee
is normal for all $\gamma,\mu\in\CC$. For the given value of $y$ those two options are the only possible ways to construct normality preserving augmentations. See the right picture in Figure \ref{figure4}.

\subsection{Illustrations belonging to Section \ref{sect3}}

\subsubsection{Normality preserving non-normal perturbation}\label{ill2}
First we illustrate Theorem \ref{Th5} by presenting an example of a line through two normal matrices that contains only two normal matrices. For this, let
\[ A=\left[\begin{array}{rr}0 & 1\\ 1 & 0 \end{array}\right], \hdrie  B = \left[\begin{array}{rr}0 & 1\\ -1 & 0\end{array}\right], \und E=B-A = \left[\begin{array}{rr}0 & 0\\ -2 & 0\end{array}\right]. \]
Thus, $E$ is a {\em non-normal normality preserving perturbation} of $A$. Hence, apart from $A$ and $B$, no other matrix of the form $A+tE$ with $t\in\RR$ is normal. Indeed,
\[ [A+tE,(A+tE)^*] = \left[\begin{array}{cc} 1-(1-2t)^2 & 0 \\ 0 & (1-2t)^2-1\end{array}\right],\]
and this matrix is only zero for $t=0$ and $t=1$. Moreover, the eigenvalues of $A+tE$ are 
\[ \sqrt{1-2t} \und -\sqrt{1-2t}, \]
whereas the corresponding sums of the eigenvalues of $\HH(A+tE)$ and $\SS(A+tE)$ equal 
\[(1-t)+it  \und -(1-t)-it.\]
Thus, for instance at $t=\half$, the eigenvalue zero of $A+tE$ is not the sum of the eigenvalues of the the Hermitian and skew-Hermitian parts of $A+tE$.

\subsubsection{Normality preserving $\theta$-Hermitian perturbations}\label{ill3}
\begin{figure}[ht]
\centering
\includegraphics[width=11cm]{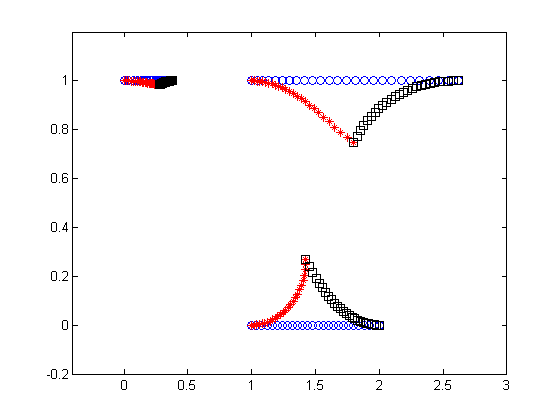}
\caption{Eigenvalue trajectories of a normality preserving perturbation, and of the same perturbation written as the sum of non-normality preserving normal perturbations.}\label{figure5}
\end{figure}

Next, we illustrate Theorem \ref{Th4}. This concerns {\em normality preserving $\theta$-Hermitian perturbations}. For this, we take $\theta=1$ and consider the matrix $A+E$, where
\[  A = \left[\begin{array}{ccc} 1 & 0 & 0 \\ 0 & i & 0 \\ 0 & 0 & 1+i\end{array}\right] \und E = u_1u_1^*+2u_2u_2^*\]
with $u_1,u_2\in\CC^3$ mutually orthonormal. Then $A+E$ is normal if and only if both $u_1$ and $u_2$ are {\em eigenvectors} of
\[ 2\Imt(A) = 2\SS(A) =  \left[\begin{array}{ccc} 0 & 0 & 0 \\ 0 & i & 0 \\ 0 & 0 & i\end{array}\right].\]
However, $A+E$, where $E=u_1u_1^*+u_2u_2^*$ with $u_1,u_2\in\CC^3$ mutually orthonormal, is normal if and only if both $u_1$ and $u_2$ are {\em linear combinations of the same two eigenvectors} $v_1$ and $v_2$ of $2\SS(A)$. Thus, with
\[ u_1 = \frac{1}{4}\left[\begin{array}{c} \sqrt{2} \\ 1 \\ 1 \end{array}\right], \hdrie u_2 = \frac{1}{4}\left[\begin{array}{c} -\sqrt{2} \\ 1 \\ 1 \end{array}\right], \hdrie\mbox{\rm where }\hdrie v_1 = \left[\begin{array}{c} 1 \\ 0 \\ 0\end{array}\right] \und  v_2 = \left[\begin{array}{c} 0 \\ \half\st \\ \half\st\end{array}\right],\] 
and $E_1=u_1u_1^*$ and $E_2=u_2u_2^*$ we have that  
\[ E_1+E_2 = \frac{1}{16}\left[\begin{array}{rcc} 2 & \st & \st \\ \st & 1 & 1 \\ \st & 1 & 1\end{array}\right]  +
\frac{1}{16}\left[\begin{array}{rcc} 2 & -\st & -\st \\ -\st & 1 & 1 \\ -\st & 1 & 1\end{array}\right]
= 
\frac{1}{16}\left[\begin{array}{ccc} 4 & 0 & 0 \\ 0 & 2 & 2 \\ 0 & 2 & 2\end{array}\right]\]
is a normality preserving rank two perturbation of $A$, written as the sum of two rank one normal perturbations $E_1$ and $E_2$ that each does not preserve normality. However, we also have that
\[ \frac{1}{4}v_1v_1^*+\half v_2v_2^* = \frac{1}{16}\left[\begin{array}{ccc} 4 & 0 & 0 \\ 0 & 0 & 0 \\ 0 & 0 & 0\end{array}\right]  +
\frac{1}{16}\left[\begin{array}{ccc} 0 & 0 & 0 \\ 0 & 2 & 2 \\ 0 & 2 & 2\end{array}\right]
=
\frac{1}{16}\left[\begin{array}{ccc} 4 & 0 & 0 \\ 0 & 2 & 2 \\ 0 & 2 & 2\end{array}\right],\]
and this expresses the same perturbation $E$ as a sum of two normality preserving normal rank-$1$ perturbations. The eigenvalues of $A+tE$ are $1+\frac{1}{4}t$ due to the term $\frac{1}{4}v_1v_1^*$, together with the eigenvalues of
\be \left[\begin{array}{cc} i & 0 \\ 0 & 1+i \end{array}\right] + \frac{1}{8}t \left[\begin{array}{cc} 1 & 1 \\ 1 & 1 \end{array}\right], \hdrie\mbox{which are }\hdrie i+\half\left( 1+\frac{t}{4} \pm \sqrt{1 + \frac{t^2}{16}}\right),\ee
which are due to the term $\frac{1}{4}v_2v_2^*$. As stated in Theorem \ref{Th4.8}, the rank-$2$ perturbation moves the eigenvalues of $A$ in the horizontal direction. Since the eigenvalues $i$ and $1+i$ of $A$ are on the same horizontal line, they can be simultaneously perturbed by a rank-$1$ perturbation. For $t\in[0,4]$ those eigenvalues are plotted by circles in Figure 5. We also computed the eigenvalues of $A+tE_1$ and of $A+4E_1+tE_2$ for $t\in[0,4]$ and indicated them by asterisks and boxes, respectively. As is visible in Figure \ref{figure5}, the eigenvalues leave the straight line before returning to the eigenvalues of the normal matrix $A+4E_1+4E_2=A+4E$.

\subsection{Illustrations belonging to Section \ref{sect5}}\label{ill4}
Finally, we will illustrate Theorem \ref{Th7}. The starting point is a $3\times 3$ matrix $A$,
\be A = \left[ \begin{array}{ccc} 5i & 0 & 0 \\ 0 & 1 & 0 \\ 0 & 0 & 2+2i\end{array}\right]\hdrie \mbox{\rm and thus, }\hdrie \HH(A)=X=\left[ \begin{array}{ccc} 0 & 0 & 0 \\ 0 & 1 & 0 \\ 0 & 0 & 2\end{array}\right]. \ee
With $\theta=1$, the eigenvalues already lie on a parabolic curve, and thus also with the trivial choice $Z=0$ augmentations $A_+$ can be constructed having eigenvalues on the same curve. More interesting is to choose a $Z\not=0$ such that $A-iZZ^*$ is normal. Since $\HH(A)$ has distinct eigenvalues, $Z$ needs to have eigenvectors of $\HH(A)$ as columns. Take for example
\be Z = \left[ \begin{array}{cc} 2 & 0 \\ 0 & 0 \\ 0 & 1\end{array}\right], \hdrie\mbox{\rm then }\hdrie ZZ^* = \left[ \begin{array}{ccc} 5 & 0 & 0 \\ 0 & 0 & 0 \\ 0 & 0 & 1\end{array}\right] \und A-iZZ^* = \left[ \begin{array}{ccc} i & 0 & 0 \\ 0 & 1 & 0 \\ 0 & 0 & 2+i\end{array}\right]. \ee
The eigenvalues of $A-iZZ^*$ lie on the curve $\mathcal{C}$ that is the image of
\be q:\RR\rightarrow \CC: \hdrie \rho \mapsto \rho + i(1-\rho)^2.\ee
Augmentations $A_+$ can now be constructed by choosing an arbitrary Hermitian $2\times 2 $ matrix $M$ and an arbitrary unitary matrix $Q$, for instance
\be M = \left[\begin{array}{cc} 1 & 1 \\ 1 & 2 \end{array}\right] \und Q = \half\sqrt{2} \left[\begin{array}{cr} 1 & 1 \\ 1 & -1 \end{array}\right],\ee
and then to form the Hermitian part of $A_+$ as
\be \HH(A_+) = \left[\begin{array}{cc} X & ZQ \\ Q^*Z^* & M \end{array}\right] = 
\left[\begin{array}{ccccc}   
 0 &  0     &    0   & \st   & \st\\
         0  &  1   &      0   &      0   &      0\\
         0  &       0   & 2   & \half\st  & -\half\st\\
    \st   &      0   & \phantom{-}\half\st  & 1   & 1\\
    \st   &      0   &-\half\st   & 1   & 2
    \end{array}\right].\ee
Finally, $A_+$ itself can be formed as $A_+ = \HH(A) + i\cdot q(\HH(A)) = \HH(A)+ i\cdot(I-\HH(A))^2$, resulting in
\be A_+ =
\left[\begin{array}{ccccc} 
        5i   & 0  &    0            &      \st      & \st + \st i\\
        0    & 1  &    0            &       0           &      0          \\
        0    & 0  &  2 + 2i         &      \hst     &-\hst - \hst i\\
   \st       & 0  & \phantom{-}\hst &    1 + 3\half i   &      1 + 2\half i\\
   \st+\st i & 0  &-\hst - \hst i   &    1 + 2\half i    &     2 + 4\half i
   \end{array}\right]. \ee
Indeed, $A_+$ is a $2$-augmentation of $A_+$. We verified that $A_+$ is normal, and computed its eigenvalues. From this example we observe that if $Z,M,Q$ are chosen real, then $\HH(A_+)$ is real symmetric and $A_+$ complex symmetric, being the sum of a real symmetric matrix and $i$ times a polynomial of this real symmetric matrix. Note that not all complex symmetric matrices are normal. In fact, the leading $4\times 4$ principal submatrix of $A_+$ in the above example is not normal, nor is the trailing $2\times 2$ principal submatrix. Thus, $A_+$ could not have been constructed using the procedure for $m=1$ twice.

\end{document}